\documentclass[12pt]{article}
\usepackage{amsfonts,amsmath,amscd}
\usepackage[a4paper,left=3cm,right=2cm,top=3cm,bottom=2cm]{geometry}
\usepackage[latin1]{inputenc}
\usepackage{graphicx}
\usepackage{wrapfig}

\newcommand{\veps}{\varepsilon}

\newcommand{\R}{\mathbb{R}}

\newtheorem{obs}{Remark}[section]

\newtheorem{teorema}{Theorem}[section]

\newcommand{\n}{\noindent}

\begin{document}

\title{A family of  warped product semi-Riemannian Einstein metrics }

\author{
\textbf{Márcio Lemes de Sousa }
\\
{\small\it ICET - CUA, Universidade Federal de Mato Grosso,}\\
{\small\it Av. Universitária nº 3.500, Pontal do Araguaia, MT,
Brazil }
\\
{\small\it e-mail:  marciolemesew@yahoo.com.br}
  \\
\textbf{Romildo Pina \footnote{ Partially supported by
CAPES-PROCAD.}}
\\
{\small\it IME, Universidade Federal de Goi\'as,}\\
{\small\it Caixa Postal 131, 74001-970, Goi\^ania, GO, Brazil }\\
{\small\it e-mail: romildo@ufg.br }
 }
\date{}

\maketitle

\thispagestyle{empty}

\markboth{abstract}{abstract}
\addcontentsline{toc}{chapter}{abstract}

\begin{abstract}
We  study warped products semi-Riemannian Einstein manifolds. We consider the case in that the base is conformal
to an $n$-dimensional pseudo-Euclidean space and invariant under
the action of an $(n-1)$-dimensional  translation group.  We provide all such solutions in the case Ricci-flat  when the base is conformal to an n-dimensional pseudo-Euclidean space, invariant under the action of an $ (n -1)$ -dimensional translation group and the fiber F is Ricci-flat. In particular, we obtain  explicit solutions, in the case vacuum, for the Einstein field equation.
 \end{abstract}

\noindent 2010 Mathematics Subject Classification: 53C21, 53C50, 53C25 \\
Key words: pseudo-Riemannian metric,  Einstein manifold, warped product

\section{Introduction and main statements}

A semi-Riemannian manifold $(M, g)$ is Einstein if there exists
a real constant $\lambda$ such that
$$Ric_{g}(X, Y) = \lambda g(X, Y)$$
for each $X, Y $ in $T_{p}M$ and each $p$ in $M$. This notion is
relevant only for $n \geq 4$. Indeed, if $n = 1$, $Ric_{g} = 0$. If
$n = 2$, then at each $p$ in $M$, we have $Ric_{g}(X, Y) =
\frac{1}{2}Kg(X, Y)$, so a 2-dimensional semi-Riemannian
manifold is Einstein if and only if it has constant sectional
or scalar curvature. If $n = 3$, then $(M, g)$ is Einstein if and
only if it has constant sectional curvature.

Over the last few years, several authors have considered the following problem:

Let $(M, g)$ be a semi-Riemannian manifold of dimension $n > 3$.
Does there exist a metric $g'$ on $M$ such that $(M, g')$ is an Einstein manifold?

According to T. Aubin \cite{Au}, to decide if a Riemannian manifold carries an Einstein
metric, will be one of the important questions in Riemannian geometry for the
next decades. Finding solutions to this problem is equivalent to solving a nonlinear system of second-order
partial differential equations. In particular, semi-Riemannian Einstein manifolds with zero Ricci curvature are special solutions, in the vacuum case $(T = 0)$, of following equation
\begin{equation}\label{eqeinst}
Ric_{g} - \frac{1}{2}K g = T,
\end{equation}
where $K$ the scalar curvature of $g$ and $T$ is a symmetric tensor of order $2$. If $g$ is the Lorentz metric on a four-dimensional manifold, this is simply the Einstein field equation. Whenever the tensor $T$ represents a physical field such as electromagnetic field perfect fluid type, pure radiation field and vacuum $(T=0)$, the above equation has been studied in several papers, most of them dealing with solutions which are invariant under some symmetry group of the equation (see \cite{KRM} for details). When the metric $g$ is conformal to the Minkowski space-time, then the solutions in the vacuum case are necessarily flat (see \cite{KRM}).

Several authors constructed new examples of Einstein manifolds. In \cite{ZIL}, Ziller constructed examples of compact
manifolds with  constant Ricci curvature.
 Chen, in \cite{CH}, constructed new examples of Einstein manifolds with odd size, and, in \cite{YAU}, Yau presented a survey on  Ricci-flat manifold.
  In \cite{Ku1}, Kühnell estudied conformal transformations between Einstein spaces and, as a consequence of
  the  obtained results, showed
  that there is no Riemannian Einstein Manifold with non-constant sectional curvature which is locally conformally flat.
   This result was extended to semi-Riemannian manifolds (see \cite{KR}).
   Accordingly to construct examples of Einstein manifolds  with non-constant sectional
   curvature, we work with manifolds that are not locally conformally
   flat. A chance to build these manifolds is to work with  warped
   product manifolds. Then considering  $(B, g_{B})$ and $(F, g_{F})$  semi-Riemannian manifolds, and let $f>0$ be a smooth function
on $B$, the warped product $M = B\times_{f}F$ is the product
manifold  $ B\times F$ furnished with the metric tensor
$$\tilde{g} = g_{B} + f^{2}g_{F},$$
$B$ is called the base of $M = B\times_{f}F$, $F$ is the fiber and
$f$ is the warping function. For example, polar coordinates
determine a warped product in the case of constant curvature
spaces, the case corresponds to $\R^{+}\times_{r}S^{n - 1} $.

There are several studies correlating warped product manifolds and
locally conformally flat manifolds, see \cite{BGV1},  \cite{BGV2},  \cite{BGV3} and their
references.

 In a series of papers, the authors studied warped product Einstein manifolds under various conditions on the curvature and symmetry,
 see \cite{QC}, \cite{CPW}and \cite{CPW1}.   Particularly,  He--Petersen--Wyle, \cite{CPW}, characterized warped product Riemannian
 Einstein metrics when the base is locally conformally flat.

It is well-known that the Einstein condition on warped geometries requires that the fibers must be
necessarily Einstein (see \cite{Be}). In this paper, initially we  give a charac--terization for warped product semi-Riemannian Einstein manifold
when the base is locally conformally flat. Using this characterization, we present new examples of semi-Riemannian manifolds with zero
Ricci curvature. More precisely, let us consider
 $(\R^n,g)$  the pseudo-Euclidean space, $n\geq 3$, with
coordinates $x=(x_1,\cdots, x_n)$ , $g_{ij}=\delta_{ij}\veps_i$ and
 $M = (\mathbb{R}^{n}, \overline{g})\times _{f}F^{m}$ a warped
product, where $\displaystyle \overline{g} =
\frac{1}{\varphi^{2}}g$,  $F$ is a semi--Riemannian
Einstein manifold  with constant Ricci curvature $\lambda_{F}$, $m\geq 1$,
$f,\varphi:\mathbb{R}^{n}\rightarrow \mathbb{R}$ are smooth
functions, where $f$ is a positive function. In Theorem 1.1, we find necessary and sufficient
conditions for the warped product metric $\widetilde{g} = \overline{g} +
f^{2}g_{F} $ to be Einstein. In Theorem 1.2, we consider $f$ and
$\varphi$ invariant under the action of an (n-1)--dimensional
translation group and  let $\xi=\sum_{i=1}^{n}\alpha_ix_i, \; \alpha_i\in
 \R$, be a basic invariant for an $(n-1)$-dimensional translation group.
We want to obtain differentiable functions $\varphi(\xi)$ and $f(\xi)$,  such that the metric
$\widetilde{g}$ is Einstein. We first  obtain necessary and sufficient conditions on $f(\xi)$ and $\varphi(\xi)$ for
the existence of $\widetilde{g}$. We show that
these conditions are different depending on the direction $\alpha=\sum_{i=1}^n \alpha_i\partial/\partial x_i$
being null (lightlike) or not. We observe that, in the null case, the metrics $\widetilde{g}$ and $g_{F}$
are necessarily Ricci-flat.

Considering $M =  (\mathbb{R}^{n}, \overline{g})\times _{f}F^{m}$ and $F$ a Ricci-flat manifold, we obtain all
the metrics $\widetilde{g} = \overline{g} + f^{2}g_{F} $ which are Ricci-flat and invariant under the action of an $(n-1)$-dimensional
translation group.
We prove that, if the direction $\alpha$ is timelike or spacelike,  the functions $f$ and $\varphi$ depend
on the dimensions $n$, $m$ and also on a finite number of parameters.
In fact, the solutions are explicitly given in Theorems 1.3 and 1.4.
 If the direction $\alpha$  is null, then there are infinitely many solutions. In fact, in this case,
for any given positive differentiable function $f(\xi)$,  the function $\varphi(\xi)$ satisfies
a linear ordinary differential equation of second order (see Theorem 1.5).  We illustrate this fact with some explicit examples. The metrics obtained in Theorems 1.3, 1.4 and 1.5 are explicit solutions, in the vacuum case, of the equation (\ref{eqeinst}). Especially considering $\tilde{g}$ Lorentz, $n = 3$ and $m = 1$ in  Theorems 1.3 and 1.5, we obtain  explicit solutions, in the case vacuum, for the Einstein field equation.

\begin{obs}
When the dimension of the fiber $F$ is $m = 1$, we consider $M =
(\R^n, \overline{g})\times _{f}\mathbb{R}$ and, in this case,
$\lambda_{F} = 0$.
\end{obs}

In what follows, we state our main results. We denote by
$\varphi_{ ,x_ix_j}$ and $f_{ ,x_ix_j}$ the second order
derivatives of $\varphi$ and $f$  with respect to $x_ix_j$.

\begin{teorema}\label{teor1}
Let $(\R^n,g)$ be a pseudo-Euclidean space, $n\geq 3$, with
coordinates $x=(x_1,\cdots, x_n)$ and $g_{ij}=\delta_{ij}\veps_i$.
Consider a warped product $M = (\mathbb{R}^{n}, \overline{g})\times _{f}F^{m}$, where $\displaystyle \overline{g} =
\frac{1}{\varphi^{2}}g$,  $F$ is a semi--Riemannian Einstein manifold
 with constant Ricci curvature $\lambda_{F}$, $m\geq 1$,
$f,\varphi:\mathbb{R}^{n}\rightarrow \mathbb{R}$ are smooth
functions, where $f$ is a positive function. Then  the warped
product metric $\widetilde{g} = \overline{g} + f^{2}g_{F} $ is Einstein
with constant Ricci curvature $\lambda$ if, and only if, the
functions
$f$ and $\varphi$ satisfy:\\

\begin{equation}\label{eqphij}
(n - 2)f\varphi_{ ,x_{i}x_{j}} - m\varphi f_{ ,x_{i}x_{j}} -
m\varphi_{ ,x_{i}}f_{ ,x_{j}} - m\varphi_{ ,x_{j}}f_{ ,x_{i}} = 0,
\qquad i\neq j,
\end{equation}

\begin{equation}\label{eqphii}
\begin{array}{rcl}
\varphi[(n-2)f\varphi_{ ,x_ix_i} - m\varphi f_{ ,x_{i}x_{i}} -
2m\varphi_{ ,x_{i}}f_{ ,x_{i}}] &+&
\varepsilon_{i}\big[f\varphi\displaystyle\sum_{k
=1}^{n}\varepsilon_{k} \varphi_{ ,x_{k}x_{k}}\\ - (n
-1)f\displaystyle\sum_{k =1}^{n}\varepsilon_{k} \varphi_{
,x_{k}}^{2}
  &+& m\varphi\displaystyle\sum_{k =1}^{n}\varepsilon_{k}
  \varphi_{ ,x_{k}}f_{ ,x_{k}}\big]\\
  &=& \varepsilon_{i}\lambda f,
\end{array}
\end{equation}
and
\begin{equation}\label{eqphjj}
-f\varphi^{2}\sum_{k =1}^{n}\varepsilon_{k} f_{ ,x_{k}x_{k}} + (n
-2)f\varphi\sum_{k =1}^{n}\varepsilon_{k} f_{ ,x_{k}}\varphi_{
,x_{k}} - (m-1)\varphi^{2}\sum_{k =1}^{n}\varepsilon_{k} f_{
,x_{k}}^{2} = \lambda f^{2} - \lambda_{F}.
\end{equation}

\end{teorema}

We want to find solutions of the system (\ref{eqphij}),
(\ref{eqphii}) and (\ref{eqphjj})
 of the form $\varphi(\xi)$ and $f(\xi)$, where   $\xi=\sum_{i=1}^{n}\alpha_ix_i, \; \alpha_i\in
 \R$.  Whenever $\sum_{i=1}^{n}\veps_i\alpha_i^2\neq 0$,  without loss of generality, we may consider
 $\sum_{i=1}^{n}\veps_i\alpha_i^2=\pm 1$.
The following theorem provides the system of ordinary differential
equations that must be satisfied by such solutions.

\vspace{.2in}

\noindent{\bf Theorem 1.2.} {\em  Let $( \R^n, g)$ be a
pseudo-Euclidean space, $n\geq 3$, with coordinates
$x=(x_1,\cdots, x_n)$ and $g_{ij}=\delta_{ij}\veps_i$.
 Consider $M = (\mathbb{R}^{n}, \overline{g})\times_{f}F^{m}$,
 where $\displaystyle \overline{g} = \frac{1}{\varphi^{2}}g$, $F^{m}$ a semi--Riemannian Einstein manifold
 with constant Ricci curvature $\lambda_{F}$ and smooth
functions $\varphi(\xi)$ and $f(\xi)$, where
$\xi=\sum_{i=1}^{n}\alpha_ix_i, \; \alpha_i\in
 \R,$ and  $\sum_{i=1}^{n}\veps_i\alpha_i^2=\veps_{i_0}$ or $\sum_{i=1}^{n}\veps_i\alpha_i^2=0$. Then $M$ is
Einstein with constant Ricci curvature $\lambda$ if, and only if,
the functions $f$ and $\varphi$ satisfy:\\ }
 \\
{\em  (i)
 \begin{equation} \label{solitonxi}
\left\{ \begin{array}{lcl}
 (n-2)f\varphi'' - m\varphi f''- 2m\varphi' f' = 0, \\
\displaystyle\sum_{k =
1}^{n}\varepsilon_{k}\alpha_{k}^{2}[f\varphi\varphi''
- (n-1)f\varphi'^{2} + m\varphi\varphi'f'] = \lambda f\\
\displaystyle\sum_{k =
1}^{n}\varepsilon_{k}\alpha_{k}^{2}[-f\varphi^{2}f'' +
(n-2)f\varphi\varphi'f' - (m-1)\varphi^{2}f'^{2}] = \lambda f^{2}
- \lambda_{F},
\end{array} \right.
 \end{equation}
  whenever }$\displaystyle\sum_{i=1}^{n}\veps_i\alpha_i^2=\veps_{i_0}, $ and  \\
{\em (ii)
\begin{equation}\label{eqxi0}
 (n-2)f\varphi''- m\varphi f''- 2m\varphi' f'=0, \quad \mbox{ \em and }\quad  \lambda = \lambda_{F}=0,
 \end{equation}
  whenever } $\displaystyle\sum_{i=1}^{n}\veps_i\alpha_i^2=0, $\\

In the following result we provide all the solutions of
(\ref{solitonxi}) when $m =1$ and $M$ is Ricci-flat.

 \vspace{.2in}
 \noindent{\bf Theorem 1.3.} {\em Let $( \R^n, g)$ be a
pseudo-Euclidean space, $n\geq 3$, with coordinates
$x=(x_1,\cdots, x_n)$ and $g_{ij}=\delta_{ij}\veps_i$.
 Consider  non-constant smooth functions $\varphi(\xi)$ and $f(\xi)$, where
$\xi=\sum_{i=1}^{n}\alpha_ix_i, \; \alpha_i\in
 \R,$ and  $\sum_{i=1}^{n}\veps_i\alpha_i^2=\veps_{i_0}$.  Then the warped product $M = (\R^n, \overline{g})\times _{f}\mathbb{R}$, with $\overline{g} = \frac{1}{\varphi^{2}}g$,
 is a Ricci-flat manifold  if, and only if,
 \begin{equation} \label{steady_xi_N}
 \left\{\begin{array}{lcl}
 \varphi(\xi)&=& \displaystyle\left[\frac{2}{(-n + 2)k_{1}\xi + k_{2}}\right]^{\frac{2}{n-2}} \\
 \mbox{} \\
 f(\xi)&=& \displaystyle\frac{2k}{(-n + 2)k _{1}\xi + k_{2}}
 \end{array} \right.
\end{equation}
{\em where $k, k_{1}$ and $k_{2}$ are constant with $k, k_{1}>0$.
These solutions are defined on the half space determined by
${\sum_{i=1}^{n}\alpha_ix_i<
\frac{k_{2}}{(n-2)k_{1}}}\cdot$}
\vskip10pt

\noindent{\bf Theorem 1.4} {\em Let $(
 \R^n, g)$ be a pseudo-Euclidean space,  $n\geq 3$, with
coordinates $x=(x_1,\cdots, x_n)$ and $g_{ij}=\delta_{ij}\veps_i$.
 Consider non-constant smooth functions $\varphi(\xi)$ and $f(\xi)$, where $\xi= \sum_{i=1}^{n}\alpha_ix_i, \; \alpha_i\in
 \R$, $\sum_{i=1}^{n}\veps_i\alpha_i^2=\veps_{i_0}$, and
  the warped product $M = (\R^n, \overline{g})\times _{f}F^{m}$, with $\overline{g} = \frac{1}{\varphi^{2}}g$
 and $F^{m}$ is a Ricci-flat manifold with $m
 \geq2$. Then $M$
 is a Ricci-flat manifold  if, and only if,}
\begin{equation} \label{steady_xi_W}
 \left\{\begin{array}{lcl}
 \varphi_{\mp}(\xi)&=& k\displaystyle\left[\mp\beta(k_{1}\xi + k_{2})\right]^{\mp\frac{\alpha}{\beta}} \\
 \mbox{} \\
 f_{\mp}(\xi)&=& \displaystyle \left[\mp\beta(k_{1}\xi + k_{2})\right]^{\mp\frac{1}{\beta}}
 \end{array} \right.
\end{equation}
{\em where $k, k_{1}, k_{2}, \alpha, \beta$ are constants with, $k,
k_{1}>0$, $\beta = \displaystyle \frac{\sqrt{m(n - 1)(m + n -
2)}}{n -1}$ and $\alpha = \displaystyle \frac{m \pm \beta}{n -
2}\cdot$} The solutions $\varphi_{-}$ and $f_{-}$ are defined on the half space determined
by ${\sum_{i=1}^{n}\alpha_ix_i < -\frac{k_{2}}{k_{1}}}$, while that $\varphi_{+}$ and $f_{+}$ are defined on the half space determined by ${\sum_{i=1}^{n}\alpha_ix_i > -\frac{k_{2}}{k_{1}}}\cdot$

 The following theorem shows that there are infinitely many
of warped products $M = (\R^n, \overline{g})\times _{f}F^{m}$
Ricci--flat, which are invariant under the action of an
$(n-1)$--dimensional group acting on $\mathbb{R}^{n}$, when $\alpha
=  \sum_{i = 1}^{n}\alpha_{i}\partial/\partial x_{i}$ is null-like vector.

\vskip10pt

\noindent{\bf Theorem 1.5} {\em Let $f(\xi)$ be any positive
differentiable function, where $\xi =  \sum_{i =
1}^{n}\alpha_{i}x_{i}$ and $ \sum_{i =
1}^{n}\varepsilon_{i}\alpha_{i}^{2} = 0$. Then there exists a
function $\varphi(\xi)$ satisfying (\ref{eqxi0}) and $M = (\R^n,
\overline{g})\times _{f}F^{m}$ is a Ricci-flat manifold. }

Before proving our main results, we give an example ilustrating
Theorem 1.5. Let $f(\xi) = ke^{A\xi}$ where $k>0$ and $A\neq 0$,
$\xi =  \sum_{i = 1}^{n}\alpha_{i}x_{i}$ and
$ \sum_{i = 1}^{n}\varepsilon_{i}\alpha_{i}^{2} = 0$.
 Solving  (\ref{eqxi0}) we get
$$\varphi(\xi) = c_{1}e^{A \left(\frac{m + \sqrt{m(n - 1)}}{n-2}\right) \xi} +
c_{2}e^{A\left(\frac{m - \sqrt{m(n - 1)}}{n-2}\right)\xi}, c_{1},
c_{2}\in\mathbb{R}.$$ It follows from Theorem 1.5 that $M = (\R^n,
\overline{g})\times _{f}F^{m}$ is a Ricci-flat manifold.
Considering $c_{1}$ and $c_{2}$ positive  constants we have that
$\varphi$ is globally defined on $\mathbb{R}^{n}$.
\begin{obs}
The condition that the warped product of an
$n$-dimensional Riemannian manifold $(M^n,g)$
with a second $m$-dimensional Einstein Riemannian manifold to be again Einstein can be
expressed as the Einstein condition

\begin{equation}\label{bak}
Ric_{f}^{m} = Ric_{g} + Hess_{g}f - \frac{1}{m} df\otimes df = \lambda g
\end{equation}

for the Bakry-Emery Ricci tensor of $g$
and the  warping function $f \in  C^{\infty}(M)$.
The   Riemannian metric $g$ and a function $f$
satisfying the
above condition is called  $ (\lambda, n+m)$ - Einstein metric.
The $ (\lambda, n+1)$ - Einstein metrics are more commonly called static metrics and such metrics have been extensively
studied for their connections
 to scalar curvature, the positive  mass theorem and general relativity. For more details see \cite{QC}.

 In this paper, in the Theorems 1.3, 1.4 and 1.5, we construct explicit examples of solutions of the
 equation (\ref{bak}) with $\lambda = 0$ in the semi-Riemannian case.
\end{obs}

\section{Proofs of the Main Results}

\n \textbf{Proof of Theorem \ref{teor1}:} Assume initially that $m>1$. It follows from
\cite{O'neil} that if $X_{1}, X_{2},\ldots, X_{n}\in{\cal
L}(\mathbb{R}^{n})$ and $Y_{1}, Y_{2}, \ldots, Y_{m}\in{\cal
L}(F)$ (${\cal L}(\mathbb{R}^{n})$ and ${\cal L}(F)$ are respectively the lift of a vector field on $\mathbb{R}^{n}$ and $F$ to $\mathbb{R}^{n}\times F$ ), then
\begin{eqnarray}
\label{ric1}
\left\{\begin{array}{lcl}
  Ric_{\widetilde{g}}(X_{i},X_{j}) &=& Ric_{\overline{g}}(X_{i},X_{j}) - \frac{m}{f}Hess_{\overline{g}}f(X_{i},X_{j}),\ \forall\ i,\ j = 1,\ldots n \\
    Ric_{\widetilde{g}}(X_{i},Y_{j}) &=& 0, \forall\ i= 1,\ldots n,\ j = 1,\ldots m  \\
    Ric_{\widetilde{g}}(Y_{i},Y_{j}) &=& Ric_{g_{F}}(Y_{i}, Y_{j}) -
  \widetilde{g}(Y_{i}, Y_{j})(\frac{\triangle_{\overline{g}}f}{f} + (m-1)\frac{\widetilde{g}(\nabla f, \nabla
  f)}{f^{2}}),\ \forall\  i,\ j = 1,\ldots m
\end{array}
\right.
\end{eqnarray}

 It is well known (see, ex. \cite{Be}) that if $\overline{g}=\frac{1}{\varphi^2}g$  , then
\[
Ric_{\overline{g}}=\frac{1}{\varphi^2}\left \{(n-2)\varphi
Hess_{g}(\varphi)+[\varphi \Delta_g
\varphi-(n-1)|\nabla_g\varphi|^2]g \right \}\,.
\]
Since $g(X_{i},X_{j}) = \varepsilon_{i}\delta_{ij}$, we have
\begin{eqnarray*}
  Ric_{\overline{g}}(X_{i}, X_{j}) &=&  \frac{1}{\varphi}\left \{(n-2)
Hess_{g}(\varphi)(X_{i}, X_{j}) \right \}\ \forall\ i\neq\ j = 1,\ldots n  \\
  Ric_{\overline{g}}(X_{i}, X_{i}) &=& \frac{1}{\varphi^2}\left \{(n-2)\varphi
Hess_{g}(\varphi)(X_{i}, X_{i})+[\varphi \Delta_g
\varphi-(n-1)|\nabla_g\varphi|^2]\varepsilon_{i} \right \}\ \forall\
i=1,\ldots n.
\end{eqnarray*}
 As $Hess_{g}(\varphi)(X_{i}, X_{j}) = \varphi_{,x_{i}x_{j}}$
, $\Delta_g \varphi = \displaystyle\sum_{k=1}^{n}\varepsilon_{k}\varphi_{,x_{k}x_{k}}$ and
$|\nabla_g\varphi|^2 = \displaystyle\sum_{k =
1}^{n}\varepsilon_{k}\varphi_{,x_{k}}^{2} $, we have
\begin{equation}\label{ric2}
\left\{ \begin{array}{ccl}
  Ric_{\overline{g}}(X_{i}, X_{j}) &=& \displaystyle\frac{(n-2)\varphi_{,x_{i}x_{j}}}{\varphi} \forall\ i\neq\ j:1\ldots\ n \\
  Ric_{\overline{g}}(X_{i}, X_{i}) &=& \displaystyle\frac{(n-2)\varphi_{,x_{i}x_{i}} +
  \varepsilon_{i}\displaystyle\sum_{k=1}^{n}\varepsilon_{k}\varphi_{,x_{k}x_{k}}}{\varphi}
  - (n - 1)\varepsilon_{i}\sum_{k=1}^{n}\frac{\varepsilon_{k}\varphi_{,x_{k}}^{2}}{\varphi^{2}}
\end{array}
\right.
\end{equation}

Recall that
\[
Hess_{\overline{g}}(f)(X_{i},X_{j})=f_{ ,x_ix_j}-\sum_k
\overline{\Gamma}_{ij}^k f_{,x_k},
\]
where $\overline{\Gamma}_{ij}^k$ are the Christoffel symbols of the metric $\overline{g}$. For $i,\ j,\ k$ distincts, we have
\[
\overline{\Gamma}_{ij}^k= 0\ \ \ \ \ \ \ \ \ \overline{\Gamma}_{ij}^i= -\frac{\varphi_{,x_{j}}}{\varphi}\ \ \ \ \ \ \ \ \ \ \overline{\Gamma}_{ii}^k= \varepsilon_{i}\varepsilon_{k}\frac{\varphi_{,x_{k}}}{\varphi}\ \ \ \ \ \ \ \ \ \overline{\Gamma}_{ii}^i= -\frac{\varphi_{,x_{j}}}{\varphi} \] therefore,

\begin{equation}\label{hes1}
\left\{
\begin{array}{ccc}
Hess_{\overline{g}}(f)(X_{i},X_{j}) &=& \displaystyle
f_{,x_ix_j}+\frac{\varphi_{,x_j}}{\varphi}f_{,x_i}
+\frac{\varphi_{,x_i}}{\varphi}f_{,x_j},\ \forall\ i\neq j = 1\ldots n\\
Hess_{\overline{g}}(f)(X_{i},X_{i}) &=& \displaystyle
f_{,x_ix_i}+2\frac{\varphi_{,x_i}}{\varphi}f_{,x_i} -
\varepsilon_{i}\sum_{k=1}^{n}\varepsilon_{k}\frac{\varphi_{,x_k}}{\varphi}f_{,x_k}
\end{array}
\right.
\end{equation}

 Substituting (\ref{ric2}) and (\ref{hes1}) in the first equation of  (\ref{ric1}), we  obtain
 \begin{equation}\label{ric3}
  Ric_{\widetilde{g}}(X_{i},X_{j}) =
 \frac{(n-2)\varphi_{,x_{i}x_{j}}}{\varphi}-\frac{m}{f}\left[f_{,x_ix_j}+\frac{\varphi_{,x_j}}{\varphi}f_{,x_i}
+\frac{\varphi_{,x_i}}{\varphi}f_{,x_j}\right],\ \forall\ i\neq j
\end{equation}

and

\begin{eqnarray}\label{ric4}
  Ric_{\widetilde{g}}(X_{i},X_{i}) &=& \frac{(n-2)\varphi_{,x_{i}x_{i}} +
  \varepsilon_{i}\displaystyle\sum_{k=1}^{n}\varepsilon_{k}\varphi_{,x_{k}x_{k}}}{\varphi}
  - (n -
  1)\varepsilon_{i}\sum_{k=1}^{n}\frac{\varepsilon_{k}\varphi_{,x_{k}}^{2}}{\varphi^{2}} \nonumber\\
  &-&  \frac{m}{f}\left[f_{,x_ix_i}+2\frac{\varphi_{,x_i}}{\varphi}f_{,x_i} -
\varepsilon_{i}\sum_{k=1}^{n}\varepsilon_{k}\frac{\varphi_{,x_k}}{\varphi}f_{,x_k}\right].
\end{eqnarray}

On the other hand,
\begin{equation}\label{gr}
\left\{
\begin{array}{ccl}
  Ric_{{g_{F}}}(Y_{i}, Y_{j}) &=& \lambda_{F}g_{F}(Y_{i}, Y_{j}) \\
  \widetilde{g}(Y_{i}, Y_{j}) &=& f^{2}g_{F}(Y_{i}, Y_{j}) \\
  \Delta_{\overline{g}}f &=& \varphi^{2}\sum_{k = 1}^{n}\varepsilon_{k}f_{,x_{k}x_{k}} - (n-2)\varphi\sum_{k = 1}^{n}\varepsilon_{k}\varphi_{,x_k}f_{,x_k} \\
  \widetilde{g}(\nabla f, \nabla f) &=& \varphi^{2}\sum_{k =
  1}^{n}\varepsilon_{k}f_{,x_{k}}^{2}
\end{array}
\right.
\end{equation}

Substituting (\ref{gr}) in the third  equation of the system (\ref{ric1}),
we have

\begin{equation}\label{ric5}
 Ric_{\widetilde{g}}(Y_{i},Y_{j}) = \gamma_{ij}g_{F}(Y_{i}, Y_{j})
\end{equation}
where, \[\gamma_{ij } =\lambda_{F} -
f\varphi^{2}\sum_{k = 1}^{n}\varepsilon_{k}f_{,x_{k}x_{k}} +
(n-2)f\varphi\sum_{k = 1}^{n}\varepsilon_{k}
\varphi_{,x_k}f_{,x_k} - (m - 1)\varphi^{2}\sum_{k =
1}^{n}\varepsilon_{k}f_{,x_{k}}^{2}  .\]

Using the equations (\ref{ric3}), (\ref{ric4}),  (\ref{ric5}) and the second
equation of (\ref{ric1}), we have $(M, \widetilde{g})$ is an Einstein manifold
if, and only if, the equations (\ref{eqphij}), (\ref{eqphii}) and
(\ref{eqphjj}) are satisfied. In the case  $m = 1$ just remember that:

\begin{eqnarray*}
  Ric_{\widetilde{g}}(X_{i},X_{j}) &=& Ric_{\overline{g}}(X_{i},X_{j}) - \frac{1}{f}Hess_{\overline{g}}f(X_{i},X_{j}),\ \forall\  i,\ j = 1,\ldots n \\
  Ric_{\widetilde{g}}(X_{i},Y) &=& 0,\ \forall\ i= 1,\ldots n\\
  Ric_{\widetilde{g}}(Y, Y) &=& - \widetilde{g}(Y, Y)\frac{\triangle_{\overline{g}}f}{f}.
\end{eqnarray*}
In this case, the equation (\ref{eqphij}) and (\ref{eqphii}) remain the same and the equation (\ref{eqphjj}) reduces to

\[ -f\varphi^{2}\sum_{k =1}^{n}\varepsilon_{k} f_{ ,x_{k}x_{k}} + (n
-2)f\varphi\sum_{k =1}^{n}\varepsilon_{k} f_{ ,x_{k}}\varphi_{
,x_{k}} = \lambda f^{2} .\]

This concludes the proof of Theorem \ref{teor1}.

\hfill $\Box$

\n \textbf{Proof of Theorem 1.2} Assume initially that $m>1$. Let  $\overline{g}=\varphi^{-2}g$ be a
conformal metric of $g$. We are assuming that $\varphi(\xi)$ and $f(\xi)$ are functions of $\xi$, where
$\xi=\sum \limits_{i=1}^n\alpha_ix_i$, $\alpha_i\in \R$ and $\sum_i\veps_i\alpha_i^2=\veps_{i_0}$ or $\sum_i\veps_i\alpha_i^2=0$.
Hence, we have
\[
\varphi_{,x_i}=\varphi'\alpha_i, \qquad
\varphi_{,x_ix_j}=\varphi''\alpha_i\alpha_j
\]
and
\[
f_{,x_i}=f'\alpha_i, \qquad f_{,x_ix_j}=f''\alpha_i\alpha_j\,.
\]
Substituting  these expressions into (\ref{eqphij}), we get
\[
(n - 2)f\varphi''\alpha_i\alpha_j - m\varphi
f''\alpha_{i}\alpha_{j} - 2m\varphi'f'\alpha_{i}\alpha_{j}= 0,
\qquad \forall\ i\neq j.
\]

If there exist $i\neq j$ such that $\alpha_i\alpha_j\neq 0$, then this equation reduces to
\begin{equation}  \label{eqphiLL}
(n - 2)f\varphi'' - m\varphi f'' - 2m\varphi'f'= 0.
\end{equation}
Similarly, considering  equation (\ref{eqphii}), we get
\[
 \alpha_i^2 \varphi\left[ (n-2)f\varphi'' - m\varphi f'' - 2m\varphi'f'  \right]
 +\veps_i\sum_k \veps_k\alpha_k^2 \left[f\varphi \varphi'' -(n-1)f(\varphi')^2 + m\varphi \varphi'f'  \right]=\veps_i\lambda f.
\]
Due to the relation between $\varphi''$ and $f''$ given in (\ref{eqphiLL}), the above equation  reduces to
\begin{equation}\label{inter}
\sum_k \veps_k\alpha_k^2 \left[f\varphi \varphi''
-(n-1)f(\varphi')^2 + m\varphi \varphi'f'  \right]=\lambda f.
\end{equation}

 Analogously, the equation (\ref{eqphjj})  reduces to
 \begin{equation}
\sum_k \veps_k\alpha_k^2 \left[-f\varphi^{2}f''
+(n-2)f\varphi\varphi' f' - (m - 1)\varphi^{2} f'^{2}  \right] = \lambda
f^{2} - \lambda_{F}.
\end{equation}

Therefore, if $\sum_k \veps_k\alpha_k^2=\veps_{i_0}$, we obtain the equations of the system (\ref{solitonxi}). If $\sum_k
\veps_k\alpha_k^2=0$, we have (\ref{eqphiLL})  satisfied and
(\ref{inter}) implies $\lambda=0$, hence $\lambda_{F} = 0$, i.e.,
(\ref{eqxi0}) holds.

If for all $i\neq\ j$, we have $\alpha_{i}\alpha_{j} = 0$, then $\xi = x_{i_{0}}$ and equation (\ref{eqphii}) is trivially satisfied for all $i\neq\ j$. Considering (\ref{eqphij}) for $i\neq\ i_{0}$, we get
\[\displaystyle\sum_{k =
1}^{n}\varepsilon_{k}\alpha_{k}^{2}[f\varphi\varphi''
- (n-1)f\varphi'^{2} + m\varphi\varphi'f'] = \lambda f\]
and hence, the second equation of (\ref{solitonxi}) is satisfied. Considering $i = i_{0}$ in (\ref{eqphij}) we get that the first equation of (\ref{solitonxi})  is satisfied.

Considering $i = i_{0}$ or $i \neq\ i_{0}$ in (\ref{eqphjj}), we get that the third equation of (\ref{solitonxi}) is  satisfied.

When $m = 1$,  the first and the second equation of the system (\ref{solitonxi}) are the same and the third equation  reduces to

\[\displaystyle\sum_{k =
1}^{n}\varepsilon_{k}\alpha_{k}^{2}[-f\varphi^{2}f'' +
(n-2)f\varphi\varphi'f'  = \lambda f^{2}.\]

This concludes the proof of Theorem 1.2.

\hfill $\Box$

\vspace{.2in}

\noindent {\bf Proof of Theorem 1.3.} We consider smooth functions $\varphi(\xi)$ and $f(\xi)$,
where $\xi=\sum_{i=1}^{n}\alpha_ix_i, \; \alpha_i\in
  \R, \; \sum_{i=1}^{n}\veps_i\alpha_i^2=\pm 1$.
It follows from Theorem 1.2 that the metric
 $\widetilde{g} = \overline{g} + f^{2}g_{F}$  is Ricci-flat if, and only if,
 $\varphi$  and $f$ satisfy
 \begin{equation} \label{steadyxi}
\left\{
 \begin{array}{rll}
 (n -2)f\varphi''- \varphi f''- 2\varphi' f' &=& 0  \\
  f\varphi\varphi''-(n-1)f(\varphi')^{2} +\varphi\varphi'f' &=&
  0\\
  \varphi f'' & = & (n - 2)\varphi'f'.
 \end{array}
\right.
 \end{equation}

 By substituting $\varphi f''$  in the first  equation of (\ref{steadyxi}), we get
 \begin{equation}\label{r1}
f\varphi'' = \frac{n}{n-2}\varphi'f'.
 \end{equation}
 Hence substituting equation (\ref{r1}) in the second equation of (\ref{steadyxi}), we have
 \[ \frac{f'}{f} = \frac{(n-2)}{2}\frac{\varphi'}{\varphi}\]
hence we get,
\begin{equation}\label{r2}
 f = k\varphi^{\frac{n-2}{2}},
\end{equation}
where $k$ is a positive constant. By deriving the equation (\ref{r2}), we obtain
\[ f' = \frac{k(n-2)}{2}\varphi^{\frac{n-4}{2}}\varphi' \]
and
\[ f'' = \frac{k(n-2)(n-4)}{4}\varphi^{\frac{n-6}{2}}(\varphi')^{2}
+ \frac{k(n-2)}{2}\varphi^{\frac{n-4}{2}}\varphi'' .\]

Substituting $f$, $f'$ e $f''$ in the first equation of (\ref{steadyxi}), we get
\[ \frac{\varphi''}{\varphi'} = \frac{n}{2}\frac{\varphi'}{\varphi}\]
hence we have
\begin{equation}\label{r3}
\varphi(\xi) = \left[\frac{2}{(-n + 2)(k_{1}\xi +
k_{2})}\right]^{\frac{2}{n-2}},
\end{equation}
where $k_{1}, k_{2}$ are constants with $k_{1} > 0$ .
By substituting (\ref{r3}) in (\ref{r2}), we obtain
\[ f(\xi) = \frac{2k}{(-n + 2)(k_{1}\xi + k_{2})}. \]
Finally, it is easily seen that $\varphi$ and $f$ satisfy the system (\ref{steadyxi}).

This concludes the proof of Theorem 1.3.

\hfill $\Box$

\noindent {\bf Proof of Theorem 1.4.}
We consider smooth functions $\varphi(\xi)$ and
$f(\xi)$, where $\xi=\sum_{i=1}^{n}\alpha_ix_i, \; \alpha_i\in
  \R, \; \sum_{i=1}^{n}\veps_i\alpha_i^2=\pm 1$.
It follows from Theorem 1.2 that the metric
 $\widetilde{g} = \overline{g} + f^{2}g_{F}$,  is Ricci-flat if, and only if,
 $\varphi$  and $f$ satisfy

 \begin{equation} \label{riccifl}
\left\{ \begin{array}{lcl}
 (n-2)f\varphi'' - m\varphi f''- 2m\varphi' f' = 0, \\
f\varphi\varphi''- (n-1)f\varphi'^{2} + m\varphi\varphi'f' = 0\\
-f\varphi^{2}f'' + (n-2)f\varphi\varphi'f' -(m-1)\varphi^{2}f'^{2}
= 0.
\end{array} \right.
 \end{equation}

The system (\ref{riccifl}) is equivalent to
\begin{equation} \label{riccifl2}
\left\{ \begin{array}{lcl}
 (n-2)\displaystyle\left(\frac{\varphi'}{\varphi}\right)' + (n-2)\displaystyle\left(\frac{\varphi'}{\varphi}\right)^{2} - m\displaystyle\left(\frac{f'}{f}\right)'
  - m\displaystyle\left(\frac{f'}{f}\right)^{2}- 2m\displaystyle\frac{\varphi'}{\varphi}\frac{f'}{f} = 0, \\
\displaystyle\left(\frac{\varphi'}{\varphi}\right)'- (n-2)\displaystyle\left(\frac{\varphi'}{\varphi}\right)^{2} +
m\displaystyle\frac{\varphi'}{\varphi}\frac{f'}{f} = 0\\
\displaystyle\left(\frac{f'}{f}\right)' -
(n-2)\displaystyle\frac{\varphi'}{\varphi}\frac{f'}{f} +
m\displaystyle\left(\frac{f'}{f}\right)^{2}= 0,
\end{array} \right.
 \end{equation}

isolating $\displaystyle\left(\frac{f'}{f}\right)'$ in the third
equation of (\ref{riccifl2}) we obtain,
\begin{equation} \label{flat}
\displaystyle\left(\frac{f'}{f}\right)' = (n -
2)\frac{\varphi'}{\varphi}\frac{f'}{f}-m\displaystyle\left(\frac{f'}{f}\right)^{2}\cdot
\end{equation}

Substituting (\ref{flat}) in the first equation of (\ref{riccifl2}) we have
\begin{equation}\label{flat2}
\displaystyle\left(\frac{\varphi'}{\varphi}\right)'
=\frac{mn}{n-2}\frac{\varphi'}{\varphi}\frac{f'}{f} -
\displaystyle\left(\frac{\varphi'}{\varphi}\right)^{2} -\frac{m(m
-1)}{n -2}\displaystyle\left(\frac{f'}{f}\right)^{2}\cdot
\end{equation}

Substituting (\ref{flat2}) in the second equation of (\ref{riccifl2})
we get
$$(n - 1)\displaystyle\left(\frac{\varphi'}{\varphi}\right)^{2} - 2 \frac{m(n-1)}{n -2}\frac{\varphi'}{\varphi}\frac{f'}{f} +
\frac{m(m -1)}{n - 2}\displaystyle\left(\frac{f'}{f}\right)^{2} = 0$$
hence we have

$$\frac{\varphi'}{\varphi} = \frac{m(n-1)\pm \sqrt{m(n - 1)(m + n
-2)}}{(n-1)(n-2)}\frac{f'}{f}\cdot$$

Making $\alpha = \displaystyle \frac{m(n-1)\pm \sqrt{m(n - 1)(m +
n -2)}}{(n-1)(n-2)}$, we have
$$\varphi = k f^{\alpha}$$
where $k>0$ is a constant. By deriving this equation we obtain \[\varphi' = k\alpha
f^{\alpha - 1}f'.\]
 Substituting $\varphi$ and $\varphi'$ in the third
equation of (\ref{riccifl2}), we get
$$\frac{f''}{f} = \left[(n - 2)\alpha - (m- 1)\right]\left(\frac{f'}{f}\right)^{2}. $$
Note that $(n-2)\alpha - (m -1) = 1 \pm\beta$, where $\beta =
\displaystyle \frac{ \sqrt{m(n - 1)(m + n -2)}}{(n-1)}$, hence we have
\begin{equation}\label{tor}
\frac{f''}{f'} = (1\pm \beta)\frac{f'}{f}\cdot
\end{equation}
Integrating the equation (\ref{tor}) we obtain:
\begin{equation}\label{fi}
f(\xi) = [\mp\beta(k_{1}\xi + k_{2})]^{\mp\frac{1}{\beta}},
\end{equation}
with $k_{1}, k_{2}$ constants and $k_{1}>0$, hence
\begin{equation}\label{fj}
\varphi(\xi) = [\mp\beta(k_{1}\xi + k_{2})]^{\mp\frac{\alpha}{\beta}}.
\end{equation}

Finally, it is easily seen that $\varphi$ and $f$ given by (\ref{fi}) and (\ref{fj}) satisfy the system (\ref{riccifl}).

This concludes the proof of Theorem 1.4.

 \hfill $\Box$

\noindent {\bf Proof of Theorem 1.5.}
Let $f(\xi)$ be any positive differentiable function invariant under translation of an given (n - 1)-dimensional translation group,
whose basic invariant $\xi=\sum_{i=1}^{n}\alpha_ix_i$, where $\alpha_i\in \R$ and  $\sum_{i=1}^{n}\veps_i\alpha_i^2 = 0$.
Then it follows from Theorem 1.2 that $M = (\R^n,
\overline{g})\times _{f}F^{m}$ is a Ricci-flat manifold, with $f$ being warping function, if, and only if,
$\lambda = \lambda_{F} = 0$ and $\varphi$ satisfies the linear ordinary differential equation (\ref{eqxi0}) determined by $f$.

 \hfill $\Box$

\end{document}